\documentclass{amsproc}

\usepackage{latexsym,amsmath,amsfonts,amscd,amssymb}
\usepackage{stmaryrd}
\usepackage{array}   
\usepackage[mathscr]{eucal} 
\usepackage{mathrsfs}
\usepackage{graphicx}
\usepackage{xypic}
\usepackage{calligra}
\usepackage[T1]{fontenc}

\DeclareFontFamily{OT1}{pzc}{}
\DeclareFontShape{OT1}{pzc}{m}{it}{<-> s * [1.10] pzcmi7t}{}
\DeclareMathAlphabet{\mathpzc}{OT1}{pzc}{m}{it}

\DeclareFontFamily{OT1}{rsfs}{}
 \DeclareFontShape{OT1}{rsfs}{n}{it}{<->rsfs10}{}
 \DeclareMathAlphabet{\curly}{OT1}{rsfs}{n}{it}

\newtheorem{theorem}{Theorem}[section]

\theoremstyle{definition}
\newtheorem{definition}[theorem]{Definition}
\newtheorem{example}[theorem]{Example}

\theoremstyle{remark}

\numberwithin{equation}{section}

\renewcommand{\leq}{\leqslant}

\renewcommand{\geq}{\geqslant}

\newcommand{\R}{\mathbb{R}}

\newcommand{\Z}{\mathbb{Z}}
\newcommand{\C}{\mathbb{C}}

\newcommand{\cM}{\mathcal{M}}

\newcommand{\calR}{\mathcal{R}}

\newcommand{\lra}{\longrightarrow}

\newcommand{\SU}{\mathrm{SU}}
\newcommand{\U}{\mathrm{U}}

\newcommand{\GL}{\mathrm{GL}}
\newcommand{\SL}{\mathrm{SL}}

\newcommand{\SSS}{\mathrm{S}}

\DeclareMathOperator{\Ad}{Ad}

\DeclareMathOperator{\rk}{rk}
\DeclareMathOperator{\rank}{rank}

\DeclareMathOperator{\Hom}{Hom}

\DeclareMathOperator{\Id}{Id}

\DeclareMathOperator{\Aut}{Aut}
\DeclareMathOperator{\Int}{Int}
\DeclareMathOperator{\Out}{Out}

\DeclareMathOperator{\reg}{reg}
\DeclareMathOperator{\op}{op}
\DeclareMathOperator{\cl}{\mathpzc{cl}}
\DeclareMathOperator{\conj}{Conj}

\renewcommand{\phi}{\varphi}

\newcommand{\liea}{\mathfrak{a}}

\newcommand{\liep}{\mathfrak{p}}
\newcommand{\liet}{\mathfrak{t}}

\newcommand{\lieg}{\mathfrak{g}}
\newcommand{\liee}{\mathfrak{e}}
\newcommand{\liek}{\mathfrak{k}}
\newcommand{\liel}{\mathfrak{l}}
\newcommand{\liez}{\mathfrak{z}}
\newcommand{\liegl}{\mathfrak{gl}}
\newcommand{\liesu}{\mathfrak{su}}
\newcommand{\lieso}{\mathfrak{so}}
\newcommand{\liesp}{\mathfrak{sp}}
\newcommand{\liesl}{\mathfrak{sl}}

\renewcommand{\phi}{\varphi}

\usepackage{tikz-cd}
\usetikzlibrary{calc}
\tikzset{curve/.style={settings={#1},to path={(\tikztostart)
			.. controls ($(\tikztostart)!\pv{pos}!(\tikztotarget)!\pv{height}!270:(\tikztotarget)$)
			and ($(\tikztostart)!1-\pv{pos}!(\tikztotarget)!\pv{height}!270:(\tikztotarget)$)
			.. (\tikztotarget)\tikztonodes}},
	settings/.code={\tikzset{quiver/.cd,#1}
		\def\pv##1{\pgfkeysvalueof{/tikz/quiver/##1}}},
	quiver/.cd,pos/.initial=0.35,height/.initial=0}
\tikzset{tail reversed/.code={\pgfsetarrowsstart{tikzcd to}}}
\tikzset{2tail/.code={\pgfsetarrowsstart{Implies[reversed]}}}
\tikzset{2tail reversed/.code={\pgfsetarrowsstart{Implies}}}
\tikzset{no body/.style={/tikz/dash pattern=on 0 off 1mm}}
\newcommand{\mtrx}[1]{\left (\begin{matrix}#1\end{matrix}\right)}

\begin{document}

\title{Vinberg pairs and Higgs bundles}

\author{Oscar Garc{\'\i}a-Prada}
\address{Instituto de Ciencias Matem\'aticas \\
  CSIC-UAM-UC3M-UCM \\ Nicol\'as Cabrera, 13--15 \\ 28049 Madrid \\ Spain}
\email{oscar.garcia-prada@icmat.es}
\thanks{
Partially supported by the Spanish Ministry of Science and
Innovation, through the ``Severo Ochoa Programme for Centres of Excellence in R\&D (CEX2019-000904-S)'' and grant
PID2019-109339GB-C31}


\date{27 August 2023}

\begin{abstract}
  We explore the role of Vinberg pairs defined by cyclic gradings of a semisimple complex Lie algebra
  in the context of Higgs bundle theory.
\end{abstract}

\maketitle

\centerline{\it Dedicated to Peter Newstead on the  occasion of his
80th birthday. }

%%%%%%%%%%%%%%%%%%%%%%%%
\section{Introduction}
%%%%%%%%%%%%%%%%%%%%%%%%
Let $G$ be a semisimple complex Lie group with Lie algebra $\lieg$ and $X$ be a compact Riemann surface of genus $g\geq 2$. A $G$-Higgs bundle over $X$ is a pair $(E,\varphi)$ consisting of a holomorphic principal $G$-bundle and a holomorphic section $\varphi$ of $E(\lieg)\otimes K_X$, where $E(\lieg)$ is the vector bundle associated to the adjoint representation of $G$ in $\lieg$, and $K_X$ is the canonical line bundle of $X$. There are suitable stability criteria for these objects and one can consider the moduli space $\cM(G)$ of polystable
$G$-Higgs bundles over $X$. Higgs bundles were introduced by Hitchin in
\cite{hitchin1987,hitchin:duke}, and their moduli spaces have been the object
of intense research for more than 35 years.

The moduli space of $G$-Higgs bundles has an extremely rich geometry. Its gauge-theoretic interpretation as  the moduli space of solutions to Hitchin's self-duality equations endows the smooth locus of the moduli space $\cM(G)$  with a hyperk\"ahler structure. The non-abelian Hodge correspondence, proved combining the   work of Hitchin  with work by Corlette \cite{corlette}, Donaldson \cite{donaldson} and Simpson \cite{simpson1}, identifies topologically $\cM(G)$ with the $G$-character variety of the fundamental group of $X$. Moreover, one has the  Hitchin map from $\cM(G)$ to an affine space obtained as a global version of  the  map $\lieg\to \lieg\sslash G$ and the Chevalley restriction theorem $\lieg\sslash G \cong \liet/W(\liet)$, where  $\liet\subset\lieg$ is a  Cartan subalgebra, and
$W(\liet)$ is the Weyl group. The generic fibres of this map are abelian varieties, defining an algebraically completely integrable system, known as the Hitchin system \cite{hitchin:duke,BNR,donagi-gaitsgory}. The Hitchin fibration plays a very important role in mirror symmetry and
Langlands duality \cite{kapustin-witten}, and has been a main ingredient in the proof of the Fundamental Lemma by Ng\^o \cite{ngo}. 

Having introduced Higgs bundles, we  introduce now Vinberg $\theta$-pairs. To do this, let $\theta$ be an automorphism of $G$ of order $m$ for an integer $m>0$. This defines an automorphism of order $m$ of $\lieg$, denoted also by $\theta$, and hence a $\Z/m\Z$-grading of $\lieg$, given by the  decomposition
\begin{equation}\label{m-grading}
\lieg =\bigoplus_{i\in \Z/m\Z} \lieg_i\;\;\mbox{with}\;\; \lieg_i=\{x\in \lieg\;\;|\;\; \theta(x)= \zeta^i x\}\;\;\mbox{and}\;\; [\lieg_i,\lieg_j]\subset \lieg_{i+j},
\end{equation}
where $\zeta$ is a primitive $m$-th root of unity.
Let $G^\theta\subset G$ be the fixed-point subgroup. This is a reductive group with Lie algebra $\lieg_0$, which acts via the adjoint representation on every $\lieg_i$ in the decomposition given in  (\ref{m-grading}). We refer to the  pairs $(G^\theta,\lieg_i)$ as Vinberg $\theta$-pairs. Sometimes we also  use this term for pairs in which $G^\theta$ is replaced by its connected component of the identity or its normalizer in $G$.   These pairs are called $\theta$-groups by Vinberg in \cite{vinberg},  as well as  Vinberg representations by  other authors \cite{levy1,levy2,levy3,panyushev2,thorne,reeder-et-al}. In \cite{wallach} Wallach calls them Vinberg pairs. Vinberg's theory is devoted to the geometric invariant theory and description of orbits  for the action of $G^\theta$ on $\lieg_i$.
The theory of Vinberg  $\theta$-pairs has been extended to  general fields of $0$ characteristic and good positive characteristic not dividing the order of $\theta$ \cite{levy2,levy3}.

Vinberg's theory is important of course in geometric invariant theory and
representation theory, and has applications to many problems. These include
the classification of  trivectors of $9$-dimensional space $\C^9$  by Elashvili--Vinberg \cite{vinberg-elashvili} using a $\Z/3\Z$-grading of $\liee_8$ given by 
$$
  \liee_8=\liesl(9,\C)\oplus \Lambda^3(\C^9)\oplus \Lambda^6(\C^9). 
$$
Certain Vinberg $\theta$-pairs have interesting connections to the arithmetic theory of elliptic curves and Jacobians (see 
  Bhargava--Gross \cite{bhargava-gross}).
There is also a connection between the Vinberg $\theta$-pair
$(\SL(9,\C)/\mu_3,\Lambda^3(\C^9))$ coming from the $\Z/3\Z$-grading of $\liee_8$ and the  moduli space of genus $2$ curves with some additional data
(see Rains--Sam \cite{rains-sam}).
Vinberg's theory has also connections to physics related to del Pezzo surfaces: mysterious duality
(see  \cite{king2}).

The object of this paper is to explore the role of Vinberg $\theta$-pairs in Higgs bundle theory.
There are two main tasks at hand. First we want to generalize certain features of the theory of $G$-Higgs bundles to Vinberg $\theta$-pairs ---note that when $m=1$, $\theta$ is the identity  and the only  Vinberg $\theta$-pair that one has in this situation is simply $(G,\lieg)$. The first step for this is to observe that, in a similar way to which we have defined $G$-Higgs bundles over $X$, we can define a $(G^\theta,\lieg_i)$-Higgs pair as a pair $(E,\varphi)$ consisting of a principal $G^\theta$-bundle $E$ over $X$ and a holomorphic section $\varphi$ of $E(\lieg_i)\otimes K_X$, where now $E(\lieg_i)$ is the vector bundle associated to the adjoint representation of $G^\theta$ in $\lieg_i$. As in the case of $G$-Higgs bundles, there are appropriate moduli spaces
$\cM(G^\theta,\lieg_i)$ of polystable  $(G^\theta,\lieg_i)$-Higgs pairs. The second step in our programme is devoted to studying how these moduli spaces appear in connection with the moduli space of $G$-Higgs bundles.

This exploration is certainly not new. There has been extensive work in the case in which $\theta$ is an involution. In fact, this was already considered by Hitchin for $G=\SL(2,\C)$ in his seminal paper \cite{hitchin1987}. A main motivation for this study has to do with the fact that the non-abelian Hodge correspondence can be extended
to the involutive case,  identifying topologically the  moduli space $\cM(G^\theta,\lieg_1)$ 
with the $G^\sigma$-character variety of the fundamental group of $X$, where $G^\sigma$ is a real form of $G$ defined by  an antiholomorphic involution  $\sigma$ of $G$ obtained by composing $\theta$ with a commuting anti-holomorphic involution $\tau$ of $G$, corresponding to a compact real form of $G$. The existence of such $\tau$ is ensured by a result of Cartan \cite{cartan}. The literature on this subject is extremely  vast (look for example  at the survey papers \cite{garcia-prada-lms,garcia-prada-handbook} for references). Of course  $(G^\theta,\lieg_0)$-Higgs pairs are just $G^\theta$-Higgs bundles  (note that, although we have defined $G$-Higgs bundles when $G$ is semisimple, the theory can be easily extended to the case in which $G$ is reductive). Both
the moduli space of $(G^\theta,\lieg_0)$-Higgs pairs
and $(G^\theta,\lieg_1)$-Higgs pairs appear inside the moduli space of $G$-Higgs bundles as fixed points of certain involutions of $\cM(G)$ (see \cite{garcia-prada-ramanan}). There are finite maps  $\cM(G^\theta,\lieg_0)\to \cM(G)$ and $\cM(G^\theta,\lieg_1)\to \cM(G)$, whose images define a hyperk\"ahler subvariety and a Lagrangian subvariety of $\cM(G)$, respectively. Here the Lagrangian condition is with respect to the natural holomorphic symplectic structure of $\cM(G)$ determined by the complex structure of $X$.

Another important feature that has also  been pursued  in the involutive case
is the generalization of the Hitchin fibration \cite{garcia-prada-peon,garcia-prada-peon-ramanan,peon,schaposnik-thesis,schaposnik}. Indeed, there is a generalization of the Hitchin map for the moduli space $\cM(G^\theta,\lieg_1)$. Here one exploits  a version  of the Chevalley restriction theorem for symmetric pairs proved by Kostant--Rallis \cite{kostant-rallis},  which establishes an isomorphism $\lieg_1\sslash G^\theta\cong \liea/W(\liea)$, where $\liea$ is a maximal abelian subalgebra of $\lieg_1$ and $W(\liea)$ is the little Weyl group. The work of Pe\'on-Nieto \cite{peon} and \cite{garcia-prada-peon} follows the  cameral approach of
Donagi--Gaitsgory \cite{donagi-gaitsgory} as well as that of Ng\^o \cite{ngo} for the case of $G$-Higgs bundles, while the work of Schaposnik \cite{schaposnik-thesis,schaposnik} for the classical groups follows the spectral curve approach initiated by Hitchin \cite{hitchin:duke}.
From both points of view one can see that, in general, the generic fibre of the Hitchin fibration is no longer abelian, except when the real form $G^\sigma$ corresponding to $\theta$ is quasi-split. This includes the split real forms (which exist for every semisimple $G$) and real forms whose Lie algebra is
isomorphic to $\liesu(n,n)$, $\liesu(n,n+1)$, $\lieso(n,n+2)$ and
$\liee_{6}(2)$. From the spectral curve perspective this is well-illustrated in \cite{hitchin-schaposnik}. Further work on the Hitchin fibration in the symmetric pair case has been done by Leslie \cite{leslie}, where the motivation here is related to a version of the Fundamental Lemma for symmetric pairs, and also more recently  by
Hameister--Morrissey \cite{hameister-morrissey}. 

Our mission in this paper is two-fold. On the one hand we review several previous results put in the larger context of Vinberg's theory of  $\theta$-pairs, and on the other, we  report on ongoing  joint work with Miguel Gonz\'alez,  where we address some general issues regarding the geometry of moduli spaces of cyclic Higgs bundles associated to general Vinberg $\theta$-pairs, and in particular the role of the Toledo invariant in this context
\cite{garcia-prada-gonzalez}.

In Section \ref{vinberg-pairs} we first recall the notion of $\Z/m\Z$-grading of a semisimple complex Lie algebra $\lieg$, and give some examples. In addition to the case of involutions, we also describe for $G=\SL(n,\C)$, the case defined by a cyclic quiver. We then define the notion of Vinberg $\theta$-pairs, associated to a $\Z/m\Z$-grading determined by an order $m$ automorphism of a semisimple complex Lie group $G$, and describe some of the main ingredientes of Vinberg's theory. A crucial one  is a Chevalley type restriction theorem for the action of $G^\theta$ on any $\lieg_i$ appearing in the $\Z/m\Z$-grading of $\lieg$, for which there is a Cartan subspace and a little Weyl group, generalizing the original Chevalley theorem for the action of $G$ on $\lieg$ and the symmetric pair case proved by Kostant--Rallis. Of course  this is the  key for the  definition below  of a  Hitchin map for the moduli spaces of
$(G^\theta,\lieg_i)$-Higgs pairs. We also review some facts studied in
\cite{garcia-prada-ramanan} about the parametrization of finite order automorphisms of $G$ up to conjugation by inner
automorphisms, and the relation to real forms of $G$ in the involutive case. We finish the section recalling the notion of quasi-split real form and proposing
a possible notion of quasi-splitness
in the higher order case.

In Section \ref{vinberg-higgs} we define Higgs bundles over a compact Riemann surface associated to Vinberg $\theta$-pairs $(G^\theta,\lieg_i)$, defined by an order $m$ automorphism $\theta$ of a semisimple complex Lie group $G$.  A particular case are $G$-Higgs bundles obtained when $m=1$, and hence  $\theta$ is the identity automorphism. After recalling the appropriate stability criteria we define the moduli spaces of
$(G^\theta,\lieg_i)$-Higgs pairs. We then review work in \cite{garcia-prada-ramanan} showing that there is a finite map from the  moduli space $\cM(G^\theta,\lieg_i)$ of
$(G^\theta,\lieg_i)$-Higgs pairs to the moduli space $\cM(G)$  of $G$-Higgs bundles, whose image is in the fixed point subvariety for an action
on $\cM(G)$ depending on $\theta$ and $\lieg_i$ of a  cyclic group of order $m$. The $G$-Higgs bundles in the image of this map are known as cyclic Higgs bundles, and have been extensively studied (see for example \cite{baraglia,collier,dai-li,garcia-prada-2007,garcia-prada-ramanan,labourie,simpson3}). We finish the section recalling the non-abelian Hodge correspondence and explaining for which
moduli spaces of $(G^\theta,\lieg_i)$-Higgs pairs one may have such correspondence.

In Section \ref{hodge-bundles} we start with  a  $\Z$-grading
$\lieg=\oplus_{i\in \Z}\lieg_i$ of $\lieg$, and consider pairs $(G_0,\lieg_i)$, where $G_0$ is the connected group corresponding to $\lieg_0$. These pairs, referred some times as Vinberg $\C^*$-pairs, 
play an essential role in the study of fixed points in $\cM(G)$ under the $\C^*$-action  --- Hodge bundles. A crucial property proved by Vinberg \cite{vinberg2} is that $\lieg_i$ is a prehomogeneous vector space for the action of $G_0$ for every $i\neq 0$. We review
the main ingredients of a theory recently introduced in \cite{BCGT}, where a Toledo invariant is associated to a Hodge bundle and a bound, similar to the one given by the Milnor--Wood inequality for the Toledo invariant of a Higgs bundle  for a real form of Hermitian type, is proved. We then consider a canonical $\Z/m\Z$-grading associated to a $\Z$-grading of $\lieg$ and extend the inequality obtained for Hodge bundles in \cite{BCGT} to the cyclic Higgs bundles for this $\Z/m\Z$-grading.

Finally, in Section \ref{hitchin-fibration}, we define the Hitchin map for the moduli space of $(G^\theta,\lieg_1)$-Higgs pairs where $\theta$ is an automorphism of $G$ of order $m>0$. We briefly  mention how one recovers for $m=1$ the original Hitchin map for the moduli space $\cM(G)$ of $G$-Higgs bundles and the Hitchin map in the symmetric pair case for $m=2$. We also comment on some of the main
difficulties to carry out the study for arbitrary $m$ \cite{garcia-prada-gonzalez}.

\noindent
    {\bf Acknowledgements:} The author wishes to thank Miguel Gonz\'alez, Nigel Hitchin and  Alastair King for very 
useful discussions. 
  
%%%%%%%%%%%%%%%%%%%%%%%%%%%%%%%%%%%%%%%%%%%%%%%%%%%%%%%%%%%%%%%%%%%%%%%%%%%%%%
\section{Graded Lie algebras and Vinberg $\theta$-pairs}\label{vinberg-pairs}
%%%%%%%%%%%%%%%%%%%%%%%%%%%%%%%%%%%%%%%%%%%%%%%%%%%%%%%%%%%%%%%%%%%%%%%%%%%%%%

%%%%%%%%%%%%%%%%%%%%%%%%%%%%%%%%%%%%%%%%%%%%%%%%%%%%%%%%%%%%%%
\subsection{$\Z/m\Z$-gradings on Lie algebras}\label{gradings}
%%%%%%%%%%%%%%%%%%%%%%%%%%%%%%%%%%%%%%%%%%%%%%%%%%%%%%%%%%%%%%

Let $m>0$ be an integer, and let $\lieg$ be a semisimple complex Lie algebra. A 
$\Z/m\Z$-{\bf grading} on $\lieg$ is a decomposition
$$
\lieg =\bigoplus_{i\in \Z/m\Z} \lieg_i\;\;\mbox{so that}\;\; [\lieg_i,\lieg_j]\subset \lieg_{i+j}.
$$
Let $\mu_m=\{z\in \C^*\;\;\mbox{such that}\;\; z^m=1\}$. Having a $\Z/m\Z$-grading on $\lieg$ is equivalent to  having a homomorphism $\tilde{\theta}: \mu_m\to \Aut(\lieg)$, where  $\lieg_i=\{x\in \lieg\;\;\mbox{such that}\;\;
\tilde{\theta}(z)x= z^i x\}$ for every $z\in \mu_m$.

One can then  choose a generator $\zeta$ of  $\mu_m$ and define an element
$\theta\in\Aut(\lieg)$ of order $m$ given by
$$
\theta(x)= \zeta^i x\;\; \mbox{for every} \;\; x\in \lieg_i.
$$
Conversely, an element $\theta\in \Aut(\lieg)$ of order $m$, defines a
$\Z/m\Z$-grading on $\lieg$ with
$$
\lieg_i=\{x\in \lieg\;\;\mbox{such that}\;\; \theta(x)= \zeta^i x\}.
$$

Let $G$ be a semisimple complex algebraic group with Lie algebra $\lieg$. We will very often assume that the $\Z/m\Z$-grading on $\lieg$ is induced from an element of $\Aut(G)$  of order $m$,   denoted also by $\theta$. This is, of course, always the case if $G$ is simply connected. Now, let $G_0\subset G$ be the connected subgroup corresponding to $\lieg_0$. The group $G_0$ is reductive
and all the subspaces $\lieg_i$ are stable under the  adjoint action of $G_0$.

%%%\begin{remark}
%%In the above we can more generally consider $\lieg$ being reductive. 
%%%\end{remark}
%%%%%%%%%%%%%%%%%%%%%%
\subsection{Examples}
%%%%%%%%%%%%%%%%%%%%%%
%\begin{example}
%{\bf Symmetric pairs}. An important
%\end{example}

\subsubsection{Symmetric pairs} \label{symmetric-pairs}
Let $m=2$. Then a $\Z/2\Z$-grading
$\lieg=\lieg_0 \oplus \lieg_1$ defines a symmetric pair. These have been extensively studied in connection with the theory of symmetric spaces and real forms of $\lieg$ and $G$.

\subsubsection{Cyclic quivers}\label{example-quiver}
Let $m\geq 2$. Let $V$ be a complex vector space equipped with a $\Z/m\Z$-grading $V=\oplus_{i\in \Z/m\Z} V_i$. Let $G=\SL(V)$. Define on $\lieg=\liesl(V)$ the  $\Z/m\Z$-grading given by
\begin{equation}\label{quiver-grading}
\lieg_i=\{A\in \liesl(V)\;\; \mbox{such that} \;\; A(V_j)\subset V_{j+i} \;\;\mbox{for every}\;\; j\in \Z/m\Z\}
\end{equation}
In this situation
$$
G_0=\SSS(\prod_{i\in  \Z/m\Z} \GL(V_i)),
$$
and the subspace $\lieg_1$ is
$$
\bigoplus_{i\in  \Z/m\Z} \Hom(V_i,V_{i+1}).
$$
We can define the {\bf quiver} $Q$ with $m$ vertices indexed by  $\Z/m\Z$ and arrows $i\mapsto i+1$ for each $i\in  \Z/m\Z$. Then $\lieg_1$ is the space of representations of $Q$ where we put $V_i$ at the vertex $i$, subject to the condition that the product of the determinants of the elements in  $\GL(V_i)$ be $1$. This can be represented by the diagram.

\[\begin{tikzcd} 	
	%% Primero establecemos los "nodos" en formato de tabla/m\Zatriz
	{V_0} &  & {V_1} &  & \dots &  & {V_{m-1}}	
	%% Ahora establecemos las flechas. 
	%% Las coordenadas de la tabla comienzan en 1-1 arriba a la izquierda.
	%% El atributo "curve" permite "doblar" la flecha para que se curve.
	\arrow["{f_0}", from=1-1, to=1-3]
	\arrow["{f_1}", from=1-3, to=1-5]
	\arrow["{f_{m-2}}", from=1-5, to=1-7]
	\arrow["{f_{m-1}}", curve={height=-24pt}, from=1-7, to=1-1].
	\end{tikzcd}\]

  For other classical groups the action of $G_0$ on $\lieg_1$ can be interpreted in terms of a cyclic quiver with some extra structure.

%%%%%%%%%%%%%%%%%%%%%%%%%%%%%%%%%%%%%%%%%%%%%%%%%%%%%%%%%%%%%%%
\subsection{Vinberg $\theta$-pairs}\label{vinberg-theta-pairs}
%%%%%%%%%%%%%%%%%%%%%%%%%%%%%%%%%%%%%%%%%%%%%%%%%%%%%%%%%%%%%%%

Let $G$ be a semisimple complex algebraic group with Lie algebra $\lieg$. Let
$\theta\in \Aut(G)$ of order $m$. As explained above, $\theta$ defines a
$\Z/m\Z$-grading $\lieg=\oplus_{i\in \Z/m\Z}\lieg_i$. Let $G_0$ be the connected subgroup of $G$ corresponding to $\lieg_0$.  As mentioned above, the group $G_0$ is reductive
and all the subspaces $\lieg_i$ are stable under the  adjoint action of $G_0$.
In  \cite{vinberg}, Vinberg studies  the invariant theory for the action of $G_0$ on   $\lieg_i$. The pairs $(G_0,\lieg_i)$ will be called
{\bf Vinberg $\theta$-pairs}. These pairs are referred by Vinberg and other authors \cite{levy1,levy2,levy3,panyushev2,thorne,reeder-et-al} as  $\theta$-groups or Vinberg representations. In \cite{wallach} Wallach calls them 
Vinberg pairs. Our terminology makes explicit the role of the automorphism $\theta$.

In  \cite{vinberg}, Vinberg observes that one can reduce such study to the action of $G_0$ on $\lieg_1$. Indeed, for any $k$ one may pass to 
a subalgebra $\lieg'$  with a new grading modulo $m'=\frac{m}{(m,k)}$, by setting
$\lieg'_j=\lieg_{jk}$ for $j\in \Z/m\Z$. We then have $\lieg_0'=\lieg_0$ and
$\lieg_1'=\lieg_k$. We can thus focus mostly on the pair $(G_0,\lieg_1)$.
Recall that if  $\liet\subset\lieg$ is a  Cartan subalgebra, and
$W(\liet)$ is the Weyl group, the {\bf Chevalley restriction theorem} establishes an isomorphism
$$
\lieg\sslash G \cong \liet/W(\liet).
$$
Similarly, if $\theta$ is an involution of $G$ and $\lieg = \lieg_0 \oplus \lieg_1$ is the Cartan decomposition
defined by $\theta$, and $W(\liea)$ is the {\bf little Weyl group} defined by a maximal abelian subalgebra $\liea\subset \lieg_1$, there is also a Chevalley restriction theorem stating
$$
\lieg_1\sslash G_0 \cong \liea/W(\liea),
$$
where $G_0\subset G$ is the connected  subgroup corresponding to $\lieg_0$.

One of the main results of Vinberg in \cite{vinberg} is a version of the Chevalley restriction theorem for a Vinberg $\theta$-pair $(G_0,\lieg_1)$ defined by a finite order element $\theta\in \Aut(G)$. A key concept is that of {\bf Cartan subspace},  a linear subspace  $\liea\subset \lieg_1$
which is abelian as a Lie algebra, consists of semisimple elements,
and is maximal with these two properties.
Any two Cartan subspaces of $\lieg_1$ are $G_0$-conjugate and any semisimple element of $\lieg_1$ is contained in a Cartan subspace. Hence the dimension of $\liea$ is an invariant of the grading, called the {\bf rank} of $\theta$. 
Given a
Cartan subspace $\liea$,   the normalizer and centralizer of $\liea$  in $G_0$, as usual, are defined respectively by 
$$
N_{G_0}(\liea)=\{g\in G_0\;|\; \Ad(g)\liea\subset \liea\},
$$
and 
$$
C_{G_0}(\liea)=\{g\in G_0\;|\; \Ad(g)x=x\;\;\;\mbox{for every}\;\; x\in \liea\}.
$$
Then the  {\bf little Weyl group} is defined by
$$
W(\liea)=N_{G_0}(\liea)/C_{G_0}(\liea).
$$
The group  $W(\liea)$ is a finite linear group  generated by semisimple transformations of $\liea$ fixing a hyperplane, and hence  $\C[\liea]^{W(\liea)}$ is a polynomial ring. As shown by Vinberg \cite{vinberg}, there is a version of the Chevalley restriction theorem, namely, the restriction of polynomial functions
  $\C[\lieg_1]\to \C[\liea]$   induces an isomorphism of invariant polynomial rings
  \begin{equation}\label{vinberg-theorem}
  \C[\lieg_1]^{G_0} \to \C[\liea]^{W(\liea)}, 
  \end{equation}
  or equivalently,
$$
\lieg_1\sslash G_0 \cong \liea/W(\liea).
$$
The fact that  $W(\liea)$ is a finite linear group  generated by complex reflections implies, by results of
Shephard--Todd \cite{shephard-todd}, and Chevalley \cite{chevalley},
that $\C[\liea]^{W(\liea)}=\C[f_1,\cdots,f_r]$ is a polynomial algebra generated by $r$ algebraically independent polynomials  $f_1,\cdots, f_r$ whose degrees $d_1,\cdots,d_r$ are determined by the grading. In particular, the product of these degrees is the order of $W(\liea)$.

In \cite{kostant} Kostant showed that the quotient map $\lieg\to \lieg\sslash G$ for the adjoint action of $G$ on $\lieg$ has a section, known as the {\bf Kostant section}. This was extended in \cite{kostant-rallis} to obtain the {\bf Kostant--Rallis section} in the symmetric pair case for the quotient map
$\lieg_1\to \lieg_1\sslash G_0$ defined by the isotropy representation. The existence of a similar section for Vinberg $\theta$-pairs for $\theta$ of higher order
was conjectured by Popov \cite{popov}. In this context, such a section is referred as a
{\bf Kostant--Weierstrass section}, as we explain below.
As outlined in \cite{vinberg-popov}, for a particularly nice action of a reductive complex algebraic group $G$ on a vector space  $V$ there exists an affine linear subvariety $W\subset V$, called a {\bf Weierstrass section} such that the restricting to $W$ induces an isomorphism $\C[V]^\C\to \C[W]$, that is, such that $W\hookrightarrow V$ is a section for he quotient morphism
$V\to V\sslash G$. This motivated the introduction of the term  {\bf Kostant--Weierstrass section} in \cite{panyushev2} for the quotient $\lieg_1\to \lieg_1\sslash G_0$ for a Vinberg $\theta$-pair $(G_0,\lieg_1)$.

The Popov conjecture was proved first  by
Panyushev in \cite{panyushev1} when $G_0$ is semisimple, and in \cite{panyushev2} when $\lieg_1$ contains a regular nilpotent element of $\lieg$. In \cite{levy1} Levy studies Vinberg $\theta$-pairs for fields of (good) positive characteristic and  extends the Kostant-Rallis results when $\theta$ is an involution, and in \cite{levy2,levy3} he deals with the existence of the Kostant--Weierstrass section in the higher order case, including the case when the field is of zero characteristic. In particular, in
\cite{levy2} he shows the existence of a  Kostant--Weierstrass section in zero and (good) positive characteristic
for the classical groups, and  in \cite{levy3}, he shows the existence of such section for $G$ of type $G_2$, $F_4$ and $D_4$.
The general case is proved in \cite{reeder-et-al} by Reeder--Levy--Yu--Gross. In \cite{reeder-et-al},  these authors give also a complete classification of positive-rank gradings and their little Weyl groups for simple algebraic groups $G$.

%%%%%%%%%%%%%%%%%%%%%%%%%%%%%%%%%%%%%%%%%%%%%%%%%%%%%%%%%%%%
\subsection{Extended Vinberg $\theta$-pairs}\label{extended}
%%%%%%%%%%%%%%%%%%%%%%%%%%%%%%%%%%%%%%%%%%%%%%%%%%%%%%%%%%%%

As above, let $G$ be a semisimple complex algebraic group, and $\theta\in \Aut(G)$ of finite order. Let $(G_0,\lieg_i)$ be a Vinberg $\theta$-pair, where
$\lieg_i$  is a subspace appearing  in the cyclic grading of $\lieg$ defined by $\theta$
. Let $G^\theta\subset G$ be the subgroup of fixed points under the action of $\theta$ on $G$. The group $G^\theta$ is reductive and has $G_0$ as the identity component. We thus have an extension
$$
1 \lra G_0\lra G^\theta \lra \Gamma^\theta \lra 1,
$$
where $\Gamma^\theta$ is the group of connected components, which is a finite group. Another relevant group is 
$$
G_\theta:=\{g\in G\;:\; \theta(g)=c(g)g,\;\mbox{with}\; c(g)\in Z(G)\},
$$
which is the  normalizer of $G^\theta$  in $G$ (see \cite{garcia-prada-ramanan})
for which there is an exact sequence
$$
1 \lra G^\theta\lra G_\theta \lra \Gamma_\theta \lra 1,
$$
where $\Gamma_\theta$ is a finite group. We can combine these two extensions to have the extension
$$
1 \lra G_0\lra G_\theta \lra \ \hat{\Gamma}_\theta \lra 1,
$$
where  $\hat{\Gamma}_\theta$ is the group of connected components of $G_\theta$.

Now, the subspaces $\lieg_i$ are  stable under the adjoint action of both $G^\theta$ and $G_\theta$. We can thus consider the pairs $(G^\theta,\lieg_i)$ and
$(G_\theta,\lieg_i)$, to which we will  refer also as Vinberg  $\theta$-pairs, using the more precise term 
{\bf extended Vinberg $\theta$-pairs} if needed. For these pairs, there are analogous result to (\ref{vinberg-theorem}), where the little Weyl group is defined replacing $G_0$ by $G^\theta$ and $G_\theta$.

%%%%%%%%%%%%%%%%%%%%%%%%%%%%%%%%%%%%%%%%%%%%%%%%%%%%%%%%%%%%%%%%%%%%
\subsection{Finite order automorphisms of $G$ and the clique map}\label{clique}
%%%%%%%%%%%%%%%%%%%%%%%%%%%%%%%%%%%%%%%%%%%%%%%%%%%%%%%%%%%%%%%%%%%%

Let $G$ be a semisimple complex algebraic group. Let $\Int(G)$ be the normal subgroup of $\Aut(G)$ consisting of inner automorphisms of $G$ and $\Out(G)=\Aut(G)/\Int(G)$. Given  $\theta,\theta'\in \Aut(G)$, we say  that $\theta \sim \theta'$ if there is an element $g\in G$ such that $\theta'=\Int_g\theta \Int_{g^{-1}}$. This defines an equivalence relation in $\Aut(G)$, so that that the quotient homomorphism $\pi:\Aut(G)\to \Out(G)$ descends to a map
\begin{equation}\label{equi-map}
\Aut(G)/\sim \to \Out(G).
\end{equation}

Let $m>0$ be a integer and let $\Aut_m(G)\subset \Aut(G)$  be the set of elements of order $m$.  If now $\theta,\theta'\in \Aut_m(G)$ and $\theta'=\Int_g\theta \Int_{g^{-1}}$, the two $\Z/m\Z$-gradings of $\lieg$ given by
$\lieg=\oplus_{i\in \Z/m\Z}\lieg_i$ and  $\lieg=\oplus_{i\in \Z/m\Z}\lieg_i'$ are related by $\lieg_i'=\Ad(g)\lieg_i$. In particular, if $G_0$ and $G_0'$ are the connected subgroups of $G$ corresponding to $\lieg_0$ and  $\lieg_0'$, respectively, then $G_0'=\Int_gG_0$ (similarly, $G^{\theta'}=\Int_gG^\theta$ and
$G_{\theta'}=\Int_gG_\theta$). 

The restriction of the map (\ref{equi-map}) to  $\Aut_m(G)$ defines a map
$$
\cl: \Aut_m(G)/\sim \to \Out_m(G)
$$
called the {\bf clique map}, where here $\Out_m(G)$ is the image of $\Aut_m(G)$ under the map  $\pi:\Aut(G)\to \Out(G)$. For an element $a\in \Out_m(G)$ we refer to the set   $\cl_m^{-1}(a)$ as  the {\bf  clique} defined by  $a$ (or  $m$-clique if we need to be more precise). To identify $\cl_m^{-1}(a)$, let  $\theta\in \Aut_m(G)$ so that $\pi(\theta)=a$ and  consider the set
$$
S_\theta:=\{s\in G\;:\; s\theta(s)\cdots \theta^{m-1}(s)=z\in Z\}.
$$
There is an  action of   $Z$, the centre of $G$,  on $S_\theta$ by multiplication and of $G$ given by
$$
s\cdot g:=g^{-1}s \theta(g)=\;\; g\in G, s\in S_\theta.
$$

In Section 2 of \cite{garcia-prada-ramanan} it is proved that the map $S_\theta\to\Aut_m(G)$ given by $s\mapsto \Int_s\theta$ defines a bijection 
$$
S_\theta/(Z\times G) \longleftrightarrow \cl^{-1}(a).  
$$ 
Moreove, there is  an interpretation in terms of non abelian Galois 
cohomology since 
$$
S_\theta/(Z\times G)= H^1(\Z/m\Z,\Ad(G)),
$$
where $\Z/m\Z$ acts on $\Ad(G)$ via $\theta$.

%%%%%%%%%%%%%%%%%%%%%%%%%%%%%%%%%%%%%%%%%%%%%%%%%%%%%%%%%%%%%%%%%%%%%%%%%%%%%%%
\subsection{Involutions and real forms of complex Lie groups}\label{real-forms}
%%%%%%%%%%%%%%%%%%%%%%%%%%%%%%%%%%%%%%%%%%%%%%%%%%%%%%%%%%%%%%%%%%%%%%%%%%%%%%%

Let $G$ be a complex Lie group. We  say that a  real Lie subgroup $G_\R\subset G$ (here we are considering the underlying real structure of $G$)
is a {\bf real form} of $G$ if  $G_\R=G^\sigma$, the fixed point set of a 
{\bf conjugation} (antiholomorphic involution)  $\sigma$ of $G$. 

Let   $\conj(G)$ be the set of conjugations of $G$.  We can define in $\conj(G)$ a similar equivalence relation to the one defined on $\Aut(G)$ in Section
\ref{clique}. 
Given  $\sigma,\sigma'\in \conj(G)$, we say  that  $\sigma\sim \sigma'$ if there is an element $g\in G$ such that $\sigma'=\Int_g\sigma \Int_{g^{-1}}$.

Let now $G$ be semisimple. A basic result due to Cartan \cite{cartan} establishes that given a compact
conjugation
$\tau\in \conj(G)$, i.e. a conjugation $\tau$ so that $G^\tau$ is a maximal compact subgroup of $G$, and a class $\alpha\in \conj(G)/\sim$, there exists a representative $\sigma$ of $\alpha$ so that $\tau\sigma=\sigma\tau$. The same applies to a class in $\Aut_2(G) /\sim$, where $\Aut_2(G)$ is the set of holomorphic involutions of $G$. We thus have that the map $\sigma\mapsto \theta:=\sigma\tau$ gives a bijection 
$$
\conj(G)/\sim \longleftrightarrow \Aut_2(G) /\sim.
$$
Combining this bijection with the map
$\cl: \Aut_2(G)/\sim \to \Out_2(G)$ (defined in Section \ref{clique}), we obtain a map
$$
\widehat{\cl}: \conj(G)/\sim \to \Out_2(G).
$$
The set $\widehat{\cl}^{-1}(1)$ consists of the equivalence classes of real
forms of {\bf Hodge type}. For example, for $G=\SL(n,\C)$, these are the real forms $\SU(p,q)$, with $p+q=n$, including the compact form
$\SU(n)$.

%%%%%%%%%%%%%%%%%%%%%%%%%%%%%%%%%%%%%%%%%%%%%%%%%%%%%%%%%%%%%%%%%%%%%
\subsection{Quasi-split $\Z/m\Z$-gradings}\label{quasi-split-gradings}
%%%%%%%%%%%%%%%%%%%%%%%%%%%%%%%%%%%%%%%%%%%%%%%%%%%%%%%%%%%%%%%%%%%%%%

In the symmetric pair situation of Example \ref{symmetric-pairs} where $m=2$ and we have $\Z/2\Z$-grading
$\lieg=\lieg_0 \oplus \lieg_1$, there is the special class of quasi-split involutions $\theta$, or equivalently quasi-split real forms of $\lieg$ or $G$. Given an involution $\theta$ of $G$ with $\Z/2\Z$-grading $\lieg=\lieg_0 \oplus \lieg_1$, as above,  let $G_0$ be the subgroup of $G$ corresponding to $\lieg_0$, and  $\liea\subset \lieg_1$ be a Cartan subspace. Let
$$
\lieg_1^{\reg}=\{x\in \lieg_1\;|\; \dim C_{G_0}(x)=\dim \liea\}
$$
be the set of {\bf regular} elements of $\lieg_1$ for the adjoint action of $G_0$.
Similarly one can define $\lieg^{\reg}$ for the adjoint action of $G$.
$$
\lieg^{\reg}=\{x\in \lieg\;|\; \dim C_{G}(x)=\dim \liet\},
$$
where $\liet$ is a Cartan subalgebra of $\lieg$.

The centralizer of $\liea$ in $\lieg_0$ is defined as
$$
C_{\lieg_0}(\liea)=\{x\in \lieg_0\;|\; [x,y]=0\;\;\;\mbox{for every}\;\; y\in \liea\}.
$$
One defines similarly $C_{\lieg_1}(\liea)$ and $C_{\lieg}(\liea)$.

There are several equivalent definitions for
$\theta$ or the pair $(G_0,\lieg_1)$ to be  quasi-split (see for example \cite{garcia-prada-peon}). The involution $\theta$ is {\bf quasi-split} if any of the following conditions is satisfied:

  (1) $C_{\lieg_0}(\liea)$ is abelian. In this case
  $C_{g}(\liea)=C_{\lieg_0}(\liea)\oplus C_{\lieg_1}(\liea)$ is a Cartan subalgebra of $\lieg$.

  (2) There exists a Borel subgroup $B\subset G$ with a $\theta$-invariant maximal torus $T\subset G$ such that $\theta(B)=B^{\op}$, where $B^{\op}$ is the opposite
  Borel subgroup.

  (3) $\lieg_1^{\reg}=\lieg_1\cap \lieg^{\reg}$.

  For $\lieg$ simple the the real forms corresponding to quasi-split involutions are the split real forms of $\lieg$ (this correspond to  $C_{\lieg_0}(\liea)=0$)
  and the real forms (up to conjugation by an element of $G$)  $\liesu(n,n)$, $\liesu(n,n+1)$, $\lieso(n,n+2)$ and $\liee_{6}(2)$. In particular, there exist one and only one quasi-split class in the clique $\widehat{\cl}^{-1}(a)$ of  every element  $a\in \Out_2(G)$.

For $m>2$ we will make the following definition.

\begin{definition}
  Let $(G_0,\lieg_1)$ be a Vinberg $\theta$-pair, and let $\liea\subset \lieg_1$ be a Cartan subspace and  $C_{\lieg_0}(\liea)$ be defined as above. Then   $(G_0,\lieg_1)$ is said to be quasi-split if
  $C_{\lieg_0}(\liea)$ is abelian. It is called split if  $C_{\lieg_0}(\liea)=0$.
\end{definition}

It would be interesting to investigate other characterizations of these conditions along the lines of the symmetric pair case.

\begin{example}
  In the case of the cyclic quivers of Example \ref{example-quiver}, let $k$ be the minimum of  $\dim V_i$. Then the quasi-split condition for $(G_0,\lieg_1)$ is equivalent to  $\dim V_i=k$ or  $\dim V_i=k+1$ for every $i$.
  When $m=2$, these correspond of course to the quasi-split real forms $\liesu(k,k)$ and $\liesu(k,k+1)$.
  The split condition for arbitrary $m$ corresponds to  $\dim V_i=1$ for every $i$,  which if $m=2$  corresponds to the real form $\liesu(1,1)$.
\end{example}

%\begin{question}

% - Compare this with the definition given by Alastair for cyclic quivers.

% - Other definitions?

% - In \cite{wallach} Wallach defines the notion of {\bf regular Vinberg pairs}. It seems that these pairs satify condition (1) of
%   the previous proposition. But is this notion equivalent to the definition of being quasi-split that we are after in te higher order case? At least we can try to test this in the  cyclic quiver case. Sae question for some notion of regularity in \cite{thorne}.
  
%\end{question}

%%%%%%%%%%%%%%%%%%%%%%%%%%%%%%%%%%%%%%%%%%%%%%%%%%%%%%%%%%%%%%%%%%%%%%%
\section{Vinberg $\theta$-pairs and Higgs bundles}\label{vinberg-higgs}
%%%%%%%%%%%%%%%%%%%%%%%%%%%%%%%%%%%%%%%%%%%%%%%%%%%%%%%%%%%%%%%%%%%%%%%

For this this section, let $X$ be a compact Riemann surface of genus $g\geq 2$ and let $K_X$ be its canonical bundle.  

%%%%%%%%%%%%%%%%%%%%%%%%%%%%%%%%%%%%%%%%%%%%
\subsection{Higgs pairs}\label{higgs-pairs}
%%%%%%%%%%%%%%%%%%%%%%%%%%%%%%%%%%%%%%%%%%%%

Let $G$ be a complex reductive Lie group with Lie algebra $\lieg$ and nondegenerate $G$-invariant bilinear form $\langle\cdot,\cdot\rangle$. Let $\rho:G\to\GL(V)$ be a holomorphic representation. 
If $E\to X$ is a principal $G$-bundle, we will denote the $V$-bundle $E\times_G V$ associated to $E$ via the representation $\rho$ by $E(V)$.

A $(G,V)$-{\bf Higgs pair} is a pair $(E,\varphi)$ where $E$ is a holomorphic principal $G$-bundle on $X$ and $\varphi$ is a holomorphic section of $E(V)\otimes K_X$. 

When $V=\lieg$ and the representation $\rho$ is the adjoint representation, a $(G,\lieg)$-Higgs pair is called a $G$-{\bf Higgs bundle}.

If $E$ is a principal $G$-bundle and $\hat G< G$ is a subgroup, then a structure group reduction of $E$ to $\hat G$ is a section $\sigma$ of the bundle $E(G/\hat G).$ Associated to such a reduction is a principal $\hat G$-subbundle $E_\sigma\subset E$ such that $E_\sigma(G)$ is canonically isomorphic to $E.$ Let $G$ be a complex reductive Lie group and $\rho:G\to\GL(V)$ be a holomorphic representation. Let $\hat G< G$ and $\hat V\subset V$ be a $\rho(\hat G)$-invariant subspace. We say that a $(G,V)$-Higgs pair $(E,\varphi)$ reduces to a $(\hat G,\hat V)$-Higgs pair, if there is a holomorphic reduction $E_{\hat G}$ of $E$ to $\hat G$ such that $\varphi\in H^0(E_{\hat G}(\hat V)\otimes K_X)\subset H^0(E_{\hat G}(V)\otimes K_X).$

To form a moduli space of Higgs pairs, we need to define suitable notions of stability. We describe this below and refer to
\cite{garcia-prada-gothen-mundet} for more details.
Fix a maximal compact subgroup $K<G$ and let $\liek$ be its Lie algebra. An element $s\in i\liek$ defines subspaces of $V$ via the representation $\rho$
$$
V_s^0=\{v\in V~|\rho(e^{ts})v= v \}\;\;\mbox{and}\;\; V_s=\{v\in V~|~\rho(e^{ts})(v)\; \text{is bounded as $t\to\infty$}\}.
$$

When $\rho:G\to\GL(\lieg)$ is the adjoint representation, $\lieg_s=\liep_s\subset\lieg$ is a parabolic subalgebra with Levi subalgebra $\lieg_s^0=\liel_s\subset\liep_s.$ The associated subgroups $L_s<P_s$ are given by
$$
L_s=\{g\in G~|~\Ad(g)s=s\}\;\;\mbox{and}\;\; P_s=\{g\in G~|~\Ad(e^{ts})(g)\ \text{is bounded as $t\to\infty$}\}.
$$
Moreover, $s$ defines a character $\chi_s:\liep_s\to\C$ by
\[\chi_s(x)=\langle s,x\rangle \text{ for $x\in\liep_s$.}\]

Let $E$ be a $G$-bundle, $s\in i\liek$ and $P_s<G$ be the associated parabolic subgroup. A reduction of structure group of $E$ to $P_s$ is a $P_s$-subbundle $E_{P_s}\subset E$, this is equivalent to a section $\sigma\in\Gamma(E(G/P_s))$ of the associated bundle. We will denote the associated $P_s$-subbundle by $E_\sigma.$ The degree of such a reduction will be defined using Chern--Weil theory. Since $P_s$ is homotopy equivalent to the maximal compact $K_s=K\cap L_s$ of $L_s$, given a reduction of structure $\sigma$ of $E$ to $P_s$, there is a further reduction $\sigma'$ of $E$ to $ K_s$ which is unique up to homotopy. Let $E_{\sigma'}\subset E$ be the resulting $K_s$ principal bundle. The curvature $F_A$ of a connection $A$ on $E_{\sigma'}$ satisfies $F_A\in\Omega^2(X,E_{\sigma'}(\liek_s)).$ Thus, evaluating the character $\chi_s$ on the curvature  we have $\chi_s(F_A)\in\Omega(X,i\R),$ and we define the degree of $\sigma$ as
$$
   \label{eq deg of parabolic reduction}\deg E(\sigma,s)=\frac{i}{2\pi}\int_X\chi_s(F_A).
 $$

   If $q$ is a rational number such that $q\cdot\chi_s$  exponentiates to a character $\tilde\chi_s:P_s\to\C^*$ and $\sigma\in\Gamma(E(G/P_s))$ is a reduction, then $E_{\sigma}\times_{\chi_s}\C^*=E_{\sigma}(\tilde\chi_s)$ is a line bundle and
   $$
   \deg E(\sigma,s)=\frac{1}{q}\deg E_\sigma(\tilde\chi_s).
$$
 When $s$ is in the centre of $\lieg$ then $P_s=G$. In this case, the degree given in the equation above is simply the degree of $E$ with respect to
  $\chi=\chi_s$,  and will be denoted $\deg_{\chi}(E)$. Again, if a multiple
  $q\cdot\chi$ exponentiates to a character $\tilde \chi:G\to\C^*$ we have 
\begin{equation}
  \label{degree-central} 
  \deg_\chi(E)=\frac{1}{q}\deg E(\tilde\chi).
\end{equation}

Let $d\rho:\lieg\to\liegl(V)$ be the differential of $\rho$ and let $\liez_\liek$ be the centre of $\liek$ and
$\liek=\liek_{ss}+\liez_{\liek}$. Define 
\[\liek_\rho=\liek_{ss}\oplus (\ker(d\rho|_{\liez_{\liek}}))^\perp.\]
We are now ready to define $\alpha$-stability notions for $\alpha\in i\liez_\liek.$

Let $\alpha\in i\liez_{\liek}$. A $(G,V)$-Higgs pair $(E,\varphi)$ is:
\begin{itemize}
  \item {\bf $\alpha$-semistable} if for any $s\in i\liek$ 
    and any holomorphic reduction $\sigma\in H^0(E(G/P_s))$ such that  $\varphi\in H^0(E_\sigma(V_s)\otimes K_X)$, we have
    $\deg E(\sigma,s)\geq \langle\alpha,s\rangle$.

\item {\bf $\alpha$-stable} if it is $\alpha$-semistable and for any $s\in i\liek_\rho$ and any holomorphic reduction $\sigma\in H^0(E(G/P_s)$ such that  $\varphi\in H^0(E_\sigma(V_s)\otimes K_X)$, we have $\deg E(\sigma,s)> \langle\alpha,s\rangle$.
 
\item {\bf $\alpha$-polystable} if it is $\alpha$-semistable and whenever
$\sigma\in H^0(E(G/P_s))$ for $s\in i\liek$ satisfy $\varphi\in H^0(E_\sigma(V_s)\otimes K_X)$ and $\deg E(\sigma,s)=\langle\alpha,s\rangle$,
 there is a further holomorphic reduction $\sigma'\in H^0(E_\sigma(P_s/L_s))$ such that $\varphi\in H^0(E_{\sigma'}(V_s^0)\otimes K_X)$.
\end{itemize}

The {\bf moduli space of $\alpha$-polystable $(G,V)$-Higgs pairs} over $X$ is defined as the set of isomorphism classes of $\alpha$-polystable $(G,V)$-Higgs pairs and will be denoted by $\cM^\alpha(G,V).$ A GIT construction of these spaces is given by Schmitt in \cite{schmitt} and by Simpson for the moduli space of $0$-polystable $G$-Higgs bundles \cite{simpson2}.

  When $\alpha=0$, we refer to $0$-stability simply as stability (similarly for (semi, poly)stability), and denote the moduli space by $\cM(G,V)$. 
The {\bf moduli space of polystable $G$-Higgs bundles} will be denoted by $\cM(G)$.

When $G$ is a classical group, Higgs pairs can be studied in terms of vector bundles. For example, if $G=\SL(n,\C)$, a $G$-Higgs bundle is equivalent to a pair
$(E,\varphi)$ where $E$ is holomorphic vector bundle of rank $n$ and trivial determinant, and $\varphi$ is a homomorphism $E\to E\otimes K_X$ with vanishing trace. These are the original objects introduced by Hitchin in \cite{hitchin1987}. Here the stability condition is the familiar slope condition
$$
\mu(F)<\mu(E),
$$
where $F$ is a nonzero proper subbundle of $E$ preserved by $\varphi$, and $\mu(E)$ is the
ratio between the degree of $E$ and  its rank.

%%%%%%%%%%%%%%%%%%%%%%%%%%%%%%%%%%%%%%%%%%%%%%%%%%%%%%%%%%%%%%%%%%%%%%%%%%%%%%%%%%%
\subsection{Vinberg $\theta$-pairs and cyclic Higgs bundles}\label{vinberg-cyclic}
%%%%%%%%%%%%%%%%%%%%%%%%%%%%%%%%%%%%%%%%%%%%%%%%%%%%%%%%%%%%%%%%%%%%%%%%%%%%%%%%%%%

Let $G$ be a semisimple complex algebraic group with Lie algebra $\lieg$ and centre $Z=Z(G)$. Let $\theta\in \Aut(G)$ be of order $m$.
As explained in Section \ref{vinberg-pairs},
$\theta$ defines a $\Z/m\Z$-grading $\lieg=\oplus_{i\in \Z/m\Z} \lieg_i$.  Let $G_0\subset G$ be the connected subgroup corresponding to $\lieg_0$. Recall also from  Section \ref{vinberg-pairs} that $G_0$ is reductive
  and all the subspaces $\lieg_i$ are stable under the  adjoint action of $G_0$.

If now $X$ is a compact Riemann surface,  following the construction given in
Section \ref{higgs-pairs}, we can consider $(G_0,\lieg_i)$-Higgs pairs over $X$ associated to the Vinberg $\theta$-pairs $(G_0,\lieg_i)$ and the moduli spaces $\cM(G_0,\lieg_i)$. In fact,  we can also consider Higgs pairs associated to the extended Vinberg $\theta$-pairs $(G^\theta,\lieg_i)$ and $(G_\theta,\lieg_i)$, for $G^\theta$ and $G_\theta$ as in Section \ref{extended} and the  moduli spaces
$\cM(G^\theta,\lieg_i)$ and  $\cM(G_\theta,\lieg_i)$ for $i\in \Z/m\Z$.

The moduli spaces associated to Vinberg pairs do appear naturally inside the moduli space $\cM(G)$ of $G$-Higgs bundles over $X$ as fixed points under various actions of a cyclic group $\Gamma$ of order $m$.
To explain this,
Let $\theta$ be a group automorphism of $G$. Let $E$ be a principal $G$-bundle and  $\theta(E):=E\times_{\theta} G$ be the principal $G$-bundle obtained by the extension of structure group defined by $\theta$.
Note that a Higgs field $\varphi\in H^0(X,E(\lieg)\otimes K_X)$ produces a Higgs 
field $\theta(\varphi)\in H^0(X,\theta(E)(\mathfrak{g})\otimes K_X)$, defining  an action of $\Aut(G)$ on $\cM(G)$.
Because of gauge equivalence, this action descends to an action of the group
$\Out(G)=\Aut(G)/\Int(G)$, where $\Int(G)$ is the normal subgroup of $\Aut(G)$ consisting of inner automorphisms of $G$.
Now, the group $\C^*$ acts on  $\cM(G)$ by multiplication of the Higgs field.

Let $\mu_m=\{z\in \C^*\;\;\mbox{such that}\;\; z^m=1\}$ and  let $\zeta\in \mu_m$ be a primitive $m$-th root of unity. Consider the
homomorphism $\mu_m\to \Aut(G)\times \C^*$ defined by $\zeta\mapsto (\theta,\zeta)$. Let $\Gamma$ be the image of this
homomorphism. If now $\theta$ is of order $m$, the group $\Gamma$ is isomorphic to  $\mu_m$ and acts on $\cM(G)$ by the rule
$$
(E,\varphi)\mapsto (\theta(E),\zeta\theta(\varphi)).
$$

The action of $\Gamma$ depends only of the image $a=\pi(\theta)$ in $\Out(G)$, since $\Int(G)$ acts trivially on $\cM(G)$. 
In \cite{garcia-prada-ramanan} it is shown that the moduli space $\cM(G^\theta,\lieg_1)$ maps into $\cM(G)$ by extending the structure group and that the image, denoted by $\widetilde\cM(G,\lieg_1)$, 
is in the fixed-point locus $\cM(G)^\Gamma$. In fact,  $\cM(G)^\Gamma$ contains
also $\widetilde{\cM}(G_{\theta'},\lieg_1')$
for $\theta'\in S_\theta$, where $S_\theta$ is defined in Section \ref{clique}, and
the subvariety $\widetilde{\cM}(G_{\theta'},\lieg_1')\subset \cM(G)$ depends 
only on  the class of
$\theta'\in S_\theta/(Z\times G)$, i.e. on the element $[\theta']$ in the clique defined by $a=\pi(\theta)$ (see Section \ref{clique}).

There is a partial converse to this result if we consider the smooth locus
$\cM_*(G)\subset \cM(G)$ consisting of stable and simple $G$-Higgs bundles. This states (see Theorem 6.3 of \cite{garcia-prada-ramanan}) that
$$
\cM_*(G)^\Gamma \subset 
\bigcup_{[\theta']\in S_\theta/(Z\times G)}
\widetilde{\cM}(G^{\theta'},\lieg_1').
$$

The elements in $\widetilde\cM(G^\theta,\lieg_1)$ are called {\bf cyclic} $G$-Higgs bundles and have been studied by several authors including
\cite{baraglia,collier,dai-li,garcia-prada-2007,garcia-prada-ramanan,labourie,simpson3}.
In particular if $G=\SL(n,\C)$, and $\theta\in \Aut(G)$ is an inner automorphism of order $m$, the Vinberg $\theta$-pair $(G^\theta,\lieg_1)$ is given by the cyclic quiver described in Example \ref{example-quiver}. Note that in this case, since $\SL(n,\C)$ is simply connected, $G^\theta$ is connected and hence coincides with $G_0$.

In many of the papers referred above  cyclic Higgs bundles are studied in
the case in which $\theta$ is an inner automorphism, and these correspond to fixed points for the action of $\mu_m$ on  $\cM(G)$ given by multiplication of the Higgs field. In \cite{garcia-prada-ramanan} this is studied also when $\theta$ is an outer automorphism, that is the cliques of non-trivial elements on $\Out_m(G)$ are considered. To be more precise,  it seems appropriate to refer  to the objects in  $\widetilde\cM(G^\theta,\lieg_1)$ as {\bf $\theta$-cyclic $G$-Higgs bundles}.

The moduli spaces $\cM(G^\theta,\lieg_i)$ for general $\lieg_i$ in the $\Z/m\Z$-grading of $\lieg$ do also show up as fixed points in $\cM(G)$. But now
we have to consider the
homomorphism $\mu_m\to \Aut(G)\times \C^*$ defined by $\zeta\mapsto (\theta,\zeta^i)$ and take as $\Gamma$ the image of this
homomorphism. Then, one has similar results to the ones described above for $\lieg_1$ (see \cite{garcia-prada-ramanan}).

Finally, the moduli spaces  $\widetilde\cM(G_\theta,\lieg_i)$ for $G_\theta$ do also appear as fixed points in $\cM(G)$, but now the action of $\Gamma$ involves the action on $\cM(G)$ of the group $H^1(X,Z)$ of isomorphism classes of $Z$-bundles on $X$ given by the analogue of tensoring by a line bundle in the vector bundle case. This is studied  in \cite{garcia-prada-ramanan}.

%%%%%%%%%%%%%%%%%%%%%%%%%%%%%%%%%%%%%%%%%%%%%%%%%%%%%%%%%%%%%%%%%
\subsection{Non-abelian Hodge correspondence and cyclic gradings}\label{section-nahc}
%%%%%%%%%%%%%%%%%%%%%%%%%%%%%%%%%%%%%%%%%%%%%%%%%%%%%%%%%%%%%%%%%

Let $G$ be a reductive Lie group (real or complex).
A {\bf representation} of $\pi_1(X)$ in
$G$ is a homomorphism $\rho\colon \pi_1(X) \to G$.
The set of all such homomorphisms,
$\Hom(\pi_1(X),G)$,  is an analytic  variety, which is algebraic
if $G$ is algebraic.
The group $G$ acts on $\Hom(\pi_1(X),G)$ by conjugation:
$$
(g \cdot \rho)(\gamma) = g \rho(\gamma) g^{-1}
$$
for $g \in G$, $\rho \in \Hom(\pi_1(X),G)$ and
$\gamma\in \pi_1(X)$. If we restrict the action to the subspace
$\Hom^+(\pi_1(X),g)$ consisting of reductive representations,
the orbit space is Hausdorff.  By a {\bf reductive representation} we mean
one that, composed with the adjoint representation in the Lie algebra
of $G$, decomposes as a sum of irreducible representations.
If $G$ is algebraic this is equivalent to the Zariski closure of the
image of $\pi_1(X)$ in $G$ being a reductive group.
The
{\bf moduli space of reductive representations} of $\pi_1(X)$ in $G$,
or $G$-{\bf character variety} of $\pi_1(X)$ is defined to be the orbit space
$$
\mathcal{R}(G) = \Hom^{+}(\pi_1(X),G) / G. 
$$
If $G$ is complex $\calR(G)$ coincides with the GIT quotient
$$
\mathcal{R}(G) = \Hom(\pi_1(X),G) \sslash G. 
$$

It has the structure of an  analytic variety 
which is algebraic if $G$ is algebraic and is real if $G$ is real or  
complex  if $G$ is complex.

Assume for the remaining part of this section that  $G$ is a semisimple complex Lie group. The  {\bf non-abelian Hodge correspondence} (\cite{hitchin1987,donaldson,corlette,simpson1}) establishes a homeomorphism
\begin{equation}\label{nahc}
\cM(G)\cong \calR(G).
\end{equation}

Let 
$\theta\in \Aut(G)$ of order $m$, and  $\lieg=\oplus_{i\in \Z/m\Z}\lieg_i$, the
$\Z/m\Z$-grading defined by $\theta$. Consider $G_0$, $G^\theta$ and $G_\theta$ as
in Section \ref{extended}.  The correspondence (\ref{nahc}) can be extended to
a reductive group $G$ by considering the subvariety $\cM_0(G)\subset \cM(G)$ consisting of $G$-Higgs bundles whose underlying $G$-bundle is topologically trivial. We thus have a homeomorphism
$$
\cM_0(G^\theta)\cong \calR(G^\theta),
$$
and similarly for $G_0$ and $G_\theta$.

If now $m=2$ and $\theta$ defines the $\Z/2\Z$-grading  $\lieg=\lieg_0\oplus \lieg_1$,  a generalization of the non-abelian Hodge correspondence (\ref{nahc}) to real forms (see e.g. \cite{garcia-prada-gothen-mundet,garcia-prada-ramanan}) gives a homeomorphism
$$
\cM(G^\theta,\lieg_1)\cong \calR(G^\sigma),
$$
where $G^\sigma$ is the real form of $G$ defined  the conjugation $\sigma=\theta\tau$, with $\tau$ a compact conjugation commuting with $\theta$ (see Section \ref{real-forms}).
The image of $\cM(G^\theta,\lieg_1)$ in $\cM(G)$ defines a Lagrangian subvariety with respect to the natural holomorphic symplectic structure of  $\cM(G)$.
If we replace $G^\theta$ by $G_0$ or $G_\theta$ we have to replace
$G^\sigma$, respectively  by the identity component of $G^\sigma$ or $G_\sigma$, the normalizer of $G^\sigma$ in $G$ (see \cite{garcia-prada-ramanan}).

%%%%%%%%%%%%%%%%%%%%%%%%%%%%%%%%%%%%%%%%%%%%%%%%%%%%%%%%%%%%%%%%%%%%%%%%%%%%%%%%%%%%
\section{$\Z$-gradings, Hodge bundles and cyclic Higgs bundles}\label{hodge-bundles}
%%%%%%%%%%%%%%%%%%%%%%%%%%%%%%%%%%%%%%%%%%%%%%%%%%%%%%%%%%%%%%%%%%%%%%%%%%%%%%%%%%%%

In this section $G$ is a semisimple complex Lie group with Lie algebra $\lieg$
and Killing form $B$. Here we follow \cite{BCGT} to which we refer for details. 

%%%%%%%%%%%%%%%%%%%%%%%%%%%%%%%%%%%%%%%%%%%%%%%%%%%%%%%%%%%%%%%%%%%%%%%%%%%%
\subsection{$\Z$-gradings and the Toledo character}\label{toledo-character}
%%%%%%%%%%%%%%%%%%%%%%%%%%%%%%%%%%%%%%%%%%%%%%%%%%%%%%%%%%%%%%%%%%%%%%%%%%%%

A $\Z$-{\bf grading} of a semisimple Lie algebra $\lieg$ is a decomposition 
\[\lieg=\bigoplus_{j\in\Z}\lieg_j\ \ \ \ \  \text{such that }\ \ \ \ \ [\lieg_i,\lieg_j]\subset\lieg_{i+j}.\]
The subalgebra $\liep=\bigoplus_{j\geq 0}\lieg_j$ is a parabolic subalgebra with Levi subalgebra $\lieg_0\subset\liep.$ There is an element $\zeta\in\lieg_0$ such that $\lieg_j=\{X\in\lieg~|~[\zeta,x]=jx\}$; the element $\zeta$ is called the {\bf grading element} of the $\Z$-grading. Having a $\Z$-grading on $\lieg$ is equivalent to having a homomorphism $\psi:\C^\ast\to \Aut(\lieg)$,
defined by
$$
\psi(z)|_{\lieg_j}=z^jI.
$$

Given a $\Z$-grading $\lieg=\bigoplus_{j\in\Z}\lieg_j$, let $G_0< G$ be the centralizer of $\zeta$; $G_0$ acts on each factor $\lieg_j$. An important fact due to  Vinberg \cite{vinberg2} is that for each $j\neq 0,$ $\lieg_j$ is a {\bf prehomogeneous vector space} for $G_0$.
This means that $\lieg_j$ (for  $j\neq 0$) has an open $G_0$-orbit, which
is then  necessarily unique and dense (see \cite{knapp,sato-kimura} for details on the theory of prehomogeneous vector spaces). The pairs
$(G_0,\lieg_i)$ are sometimes also referred as Vinberg pairs. To distinguish them from the Vinberg $\theta$-pairs for $\theta$ an automorphism of finite order, to the prehomogeneous vector spaces $(G_0,\lieg_i)$ obtained from a $\Z$-grading of $\lieg$ we will call them {\bf Vinberg $\C^*$-pairs}.

 Without loss of generality, we can consider the  Vinberg $\C^*$-pair $(G_0,\lieg_1)$. Let $\Omega\subset\lieg_1$ be the open $G_0$-orbit. Since $\lieg_0$
 is the centralizer of $\zeta,$ $B(\zeta,-):\lieg_0\to\C$ defines a character. 
  The {\bf Toledo character} $\chi_T:\lieg_0\to\C$ is defined by  
  \[\chi_T(x)=B(\zeta,x)B^*(\gamma,\gamma)~,\]
  where $B^*$ is the bilinear form induced by $B$ on the dual of $\lieg$, and $\gamma$ is the longest root such that $\lieg_\gamma\subset\lieg_1.$
  The normalization factor $B^*(\gamma,\gamma)$ guarantees that the Toledo character is independent of the choice of invariant bilinear form $B.$   Let $e\in\lieg_1$ and $\{f,h,e\}$ be an $\liesl_2$-triple with $h\in\lieg_0.$ We define the {\bf Toledo rank} of $e$ by 
  \[\rk_T(e)=\frac{1}{2}\chi_T(h),\]
  and the {\bf Toledo rank} of $(G_0,\lieg_1)$ by 
  \[\rk_T(G_0,\lieg_1)=\rk_T(e) \ \ \text{for $e\in\Omega$.}\]

  Given a $\Z$-grading $\lieg=\bigoplus_{j\in\Z}\lieg_j$, there is an involution $\theta:\lieg\to\lieg$ 
  defined by
  \[\theta|_{\lieg_j}=(-1)^j\Id.\]
  This corresponds to a real form of $\lieg$ of Hodge type (and of $G$ if there is a lift of $\theta$ to $G$).

%%%%%%%%%%%%%%%%%%%%%%%%%%%%%%%%%%%%%%%%%%%%%%%%%%%
\subsection{Hodge bundles and the Toledo invariant}
%%%%%%%%%%%%%%%%%%%%%%%%%%%%%%%%%%%%%%%%%%%%%%%%%%%
Let
$\lieg=\bigoplus_{j\in\Z}\lieg_j$ be a $\Z$-grading 
with grading element $\zeta\in\lieg_0$, and  let $G_0<G$ be the centralizer of $\zeta$. 
Note that $\exp(\lambda\zeta)$ is in the centre of $G_0$ and $\Ad(\exp(\lambda\zeta))$ acts on each $\lieg_j$ by
$\exp(j\lambda)\cdot\Id.$ Let $X$ be a compact Riemann surface of genus $g\geq 2$.
Let $(E_{G_0},\varphi)$ be a $(G_0,\lieg_k)$-Higgs pair over $X$. Note that the central element $\exp(\frac{\lambda}{k}\zeta)\in G_0$ defines a holomorphic automorphism of $E_{G_0}$ and acts on $\varphi$ by multiplication by $\exp(\lambda)$. As a result we have an isomorphism of $(G_0,\lieg_k)$-Higgs pairs
 \begin{equation}
   \label{eq iso of G0gk pairs}(E_{G_0},\varphi)\cong(E_{G_0},\lambda\varphi) \text{~ \ for all $\lambda\in\C^*$}.
 \end{equation}
Extending the structure group defines a $G$-Higgs bundle $(E_G,\varphi)$ 
\[(E_G,\varphi)=(E_{G_0}(G),\varphi)\]
since $E_{G_0}(\lieg_k)\subset E_{G_0}(\lieg)\cong E_{G}(\lieg)$. Moreover, $(E_G,\varphi)\cong (E_{G},\lambda\varphi)$ for all $\lambda\in\C^*.$ 

A $G$-Higgs bundle $(E,\varphi)$ is called a {\bf Hodge bundle of type} $(G_0,\lieg_k)$ if it reduces to a $(G_0,\lieg_k)$-Higgs pair. There is a natural $\C^*$-action on the moduli spaces of polystable $G$-Higgs bundles given by $\lambda\cdot(E,\varphi)=(E,\lambda\varphi).$ The $\C^*$-action is trivial on the moduli space of $(G_0,\lieg_k)$-Higgs pairs by \eqref{eq iso of G0gk pairs}. 
As a result, polystable Higgs bundles which are Hodge bundles define fixed points of the $\C^*$-action on the moduli space of Higgs bundles. Simpson proved the converse \cite{simpson1,simpson2}. Namely, all $\C^*$-fixed points in the moduli space of Higgs bundles are Hodge bundles for some $\Z$-grading. 
In fact, if a Higgs bundle is a Hodge bundle, then polystability of the Higgs bundle is equivalent to polystability of the associated pair (see \cite{BCGT}). In particular, there is a well defined map of moduli spaces
 \[\cM(G_0,\lieg_k)\to\cM(G)\]
 whose image consists of $\C^*$-fixed points, and every $\C^*$-fixed point in $\cM(G)$ is in the image of such a map for some
 $(G_0,\lieg_k)$ of some $\Z$-grading. Via the non-abelian Hodge correspondence, Hodge bundles correspond to {\bf variations of Hodge structure} (see \cite{BCGT}).

  Let $(E,\varphi)$ be a $(G_0,\lieg_1)$-Higgs pair and $\chi_T:\lieg_0\to\C$ be the Toledo character defined in Section
 \ref{toledo-character}. Then the {\bf Toledo invariant} $\tau(E,\varphi)$ is defined to by 
\begin{equation}\label{toledo-invariant}
  \tau(E,\varphi)=\deg_{\chi_T}(E),
\end{equation}
  where $\deg_{\chi_T}(E)$ is given by (\ref{degree-central}).
Recall the definition  of the Toledo rank $\rk_{\chi_T}(v)$ of a point $v\in\lieg_1$ given in Section \ref{toledo-character}. We define the {\bf Toledo rank} $\rk_T(\varphi)$ of a $(G_0,g_1)$-Higgs pair $(E,\varphi)$ to be 
 \[\rk_T(\varphi)=\rk_T(\varphi(x))\ \text{ for a generic $x\in X$}.\]
A main result proved in \cite{BCGT} is that if $(E,\varphi)$ is a semistable, then one has the inequality 
\begin{equation}\label{am-inequality}
 \tau(E,\varphi)\geq -\rk_T(\varphi)(2g-2),  
\end{equation}
referred in \cite{BCGT} as {\bf Arakelov--Milnor inequality}, since it generalizes the Milnor--Wood inequality for the Toledo invariant of Higgs bundles for real forms of Hermitian type \cite{BGR} (as we will see below), and is analogous to the Arakelov inequalities of classical variations of Hodge structure.

%%%%%%%%%%%%%%%%%%%%%%%%%%%%%%%%%%%%%%%%%%%%%%%%%%%%%%%%%%%%%%%%%%%%%%%%
\subsection{Hodge bundles and cyclic Higgs bundles}\label{hodge-cyclic}
%%%%%%%%%%%%%%%%%%%%%%%%%%%%%%%%%%%%%%%%%%%%%%%%%%%%%%%%%%%%%%%%%%%%%%%%
Let 
\begin{equation}\label{Z-grading}
\lieg=\bigoplus_{j\in\Z}\lieg_j
\end{equation}
be a $\Z$-grading with grading element $\zeta\in\lieg_0$, and  let $G_0<G$ be the centralizer of $\zeta$.
The $\Z$-grading (\ref{Z-grading})
defines a $\Z/k$-grading of $\lieg$ for any integer $k>0$, given by composing the inclusion homomorphism
$\mu_k\subset \C^\ast$ of the group $\mu_k$ of $k$-th roots of unity with the
homomorphism $\C^\ast \to \Aut(\lieg)$ corresponding to the $\Z$-grading. In particular, let $m\geq 2$ be the smallest  integer in (\ref{Z-grading}) for which $\lieg_j=0$ for every $|j|\geq m$. Then the $\Z/m\Z$-grading defined by (\ref{Z-grading}) is given by
\begin{equation}\label{Z-m-grading}
\lieg=\bigoplus_{i\in\Z/m\Z}\overline{\lieg}_i,
\end{equation}
with
$$
\overline{\lieg}_i=\lieg_i\oplus \lieg_{i-m}\;\;\mbox{for}\;\; i\in \{0,\cdots, m-1\}.
$$

We will assume that the automorphism of $\lieg$ of order $m$ defining the $\Z/m\Z$-grading (\ref{Z-m-grading}) lifts to an automorphism of $G$, and that $G^\theta=G_0$.
%\begin{question}
%What are the conditions for  $G^\theta=G_0$? Is this true, for example,  if $G$ is simply connected? 
%\end{question}

Let $X$ be a compact Riemann surface. We want to study $(G_0,\overline{\lieg}_i)$-Higgs pairs over $X$, which as explained in Section \ref{vinberg-cyclic}, correspond to $\theta$-cyclic $G$-Higgs bundles. Without loss of generality we will consider
$(G_0,\overline{\lieg}_1)$-Higgs pairs. Now, let $(E,\varphi)$ be a
$(G_0,\overline{\lieg}_1)$-Higgs pair over $X$. Notice
that
$$
E(\overline{\lieg}_1)=E(\lieg_1)\oplus E(\lieg_{1-m}),
$$
and hence, according to this, we can decompose $\varphi=\varphi^++\varphi^-$ with
$$
\varphi^+\in H^0(X,E(\lieg_1)\otimes K_X)\;\;\mbox{and}\;\;
\varphi^-\in H^0(X,E(\lieg_{1-m})\otimes K_X).
$$
We can now consider the Toledo invariant $\tau^+=\tau(E,\varphi^+)$
as  defined by  (\ref{toledo-invariant}) using the Toledo character $\chi_T^+$ of $(G_0,\lieg_1)$ determined by the grading element
$\zeta$.
As mentioned above, if $\varphi^-=0$ one has that the semistability of  $(E,\varphi)$ implies the inequality
(\ref{am-inequality}). But actually one has the following stronger result.

\begin{theorem}\label{cyclic-am}
  Let $(E,\varphi^++\varphi^-)$ be a $(G_0,\overline{\lieg}_1)$-Higgs pair over $X$ with $\varphi^+\neq 0\neq \varphi^-$.
  Then if  $(E,\varphi^++\varphi^-)$ is semistable one has
  $$
\tau^+\geq -\rk_T(\varphi^+)(2g-2). 
$$
\end{theorem}

The proof of Theorem \ref{cyclic-am} is given in \cite{garcia-prada-gonzalez} following \cite{BGR,BCGT}.

We  can also consider the Toledo invariant $\tau^-=\tau(E,\varphi^-)$
as  defined by  (\ref{toledo-invariant}) using the Toledo character $\chi_T^-$ of $(G_0,\lieg_{1-m})$. Notice that the grading element defining the Toledo character $\chi_T^-$ is $\zeta^-=\frac{\zeta}{1-m}$. As in the previous case, if  $\varphi^+=0$ one has
\begin{equation}\label{am-oposite}
\tau^-\geq -\rk_T(\varphi^-)(2g-2). 
\end{equation}
But now, in general we can not guarantee this inequality if
$(E,\varphi^++\varphi^-)$ is semistable with $\varphi^+\neq 0\neq \varphi^-$. This happens however under some conditions, in particular when $m=2$. In this case, the roles of $\lieg_1$ and
$\lieg_{1-m}=\lieg_{-1}$ are symmetric and we can apply Theorem \ref{cyclic-am}, reversing the roles of $\varphi^+$ and $\varphi^-$ to obtain (\ref{am-oposite}) when $\varphi^+\neq 0\neq \varphi^-$.

In fact, when $m=2$, $\tau:=\tau^+=-\tau^-$ and the semistability of $(E,\varphi^+,\varphi^-)$ implies

\begin{equation}\label{mw-inequality}
- \rk_T(\varphi^+)(2g-2)\leq \tau\leq  \rk_T(\varphi^-)(2g-2).
\end{equation}

When $\lieg$ admits a a $\Z$-grading of the form $\lieg=\lieg_{-1}\oplus \lieg_0\oplus \lieg_1$,
the associated  $\Z/2\Z$-grading $\lieg=\overline{\lieg}_0\oplus \overline{\lieg}_1$, with $\overline{\lieg}_0=\lieg_0$ and
$\overline{\lieg}_1=\lieg_{-1}\oplus \lieg_1$, corresponds to a real form  of Hermitian type. If $\lieg$ is simple a real form of Hermitian type is one of the classical real Lie algebras $\liesu(p,q),\liesp(2n,\R),\lieso^*(2n),\lieso(2,n)$ or
$\liee_6(-14),\liee_7(-25)$ in the exceptional case. In this situation, if the involution of $\lieg$ can be lifted to an involution $\theta$ of  $G$ (which is necessarily inner) the moduli space $\cM(G^\theta,\overline{\lieg}_1)$ maps to the fixed point subvariety  of the
involution of $\cM(G)$ given by $(E,\varphi)\mapsto (E,-\varphi)$ which  is homeomorphic, via the non-abelian Hodge correspondence, to the character variety $\calR(G^\sigma)$, where $\sigma$ is an antiholomorphic involution of $G$ corresponding to $\theta$, and the symmetric space defined by the quotient  $G^\sigma$ by its maximal compact subgroup is a Hermitian symmetric space of the non-compact type.  In this case, one has the  Milnor--Wood inequality  (\ref{mw-inequality}) for a semistable
$(G^\theta,\overline{\lieg}_1)$-Higgs pair $(E,\varphi^+,\varphi^-)$ as proved in \cite{BGR}. This inequality leads to the classical
Milnor--Wood inequality
$$
|\tau|\leq 2r(g-1),
$$
where $r$ is the rank of the Hermitian symmetric space. Here  $\tau$ coincides with  the original invariant defined by Toledo for a representation of $\pi_1(X)$ in $G^\sigma$ for which the above bound is proved in general by Burger--Iozzi--Wienhard
\cite{BIW}.

%%%%%%%%%%%%%%%%%%%%%%%%%%%%%%%%%%%%%%%%%%%%%%%%%%%%%%%%
\subsection{Chains  and cyclic $\SL(n,\C)$-Higgs bundles}
%%%%%%%%%%%%%%%%%%%%%%%%%%%%%%%%%%%%%%%%%%%%%%%%%%%%%%%%%

A particular example of a $\Z/m\Z$-grading of a Lie algebra $\lieg$ arising from a $\Z$-grading of $\lieg$ in the way described in Section \ref{hodge-cyclic} is given by the cyclic quiver given in  Example \ref{example-quiver}.  In this example $V$ is a complex vector space equipped with a $\Z/m\Z$-grading $V=\oplus_{i\in \Z/m\Z} V_i$ and  $G=\SL(V)$. On $\lieg=\liesl(V)$ there is a $\Z$-grading for which 
$$
G_0=\SSS(\prod_{i\in  \Z/m\Z} \GL(V_i)),
$$
and 
$$
\lieg_1=\bigoplus_{i=1}^{m-1} \Hom(V_{i-1},V_{i})\;\; \mbox{and}\;\; \lieg_{1-m}=\Hom(V_{m-1},V_0).
$$
From here we obtain the $\Z/m\Z$-grading of $\lieg$ described in Example \ref{example-quiver}, for which
$$
\overline{\lieg}_1=\lieg_1\oplus \lieg_{1-m}.
$$
Let $X$ be a compact Riemann surface of genus $g\geq 2$ with canonical bundle $K_X$. A  $(G_0,\overline{\lieg}_1)$-Higgs pair over $X$ is equivalent to the
$K_X$-twisted quiver bundle over $X$

\begin{equation}\label{quiver-bundle}
\begin{tikzcd} 	
	%% Primero establecemos los "nodos" en formato de tabla/matriz
	{E_0} &  & {E_1} &  & \dots &  & {E_{m-1}}	
	%% Ahora establecemos las flechas. 
	%% Las coordenadas de la tabla comienzan en 1-1 arriba a la izquierda.
	%% El atributo "curve" permite "doblar" la flecha para que se curve.
	\arrow["{\varphi_0}", from=1-1, to=1-3]
	\arrow["{\varphi_1}", from=1-3, to=1-5]
	\arrow["{\varphi_{m-2}}", from=1-5, to=1-7]
	\arrow["{\varphi_{m-1}}", curve={height=-24pt}, from=1-7, to=1-1],
\end{tikzcd}
\end{equation}
with $E_i$ holomorphic vector bundles over $X$ of rank $n_i=\dim V_i$ and
homomorphisms 
$\varphi_i:E_i\to E_{i+1}\otimes K_X$ for $i=0,\cdots,m-1$. To such an object we can associate a $G$-Higgs bundle
$(E,\varphi)$ with
$$
E=E_0\oplus\cdots \oplus  E_{m-1}\;\;\; \mbox{and}\;\;\;
\varphi =
\mtrx{0& 0&0&\hdots&\varphi_{m-1}\\
  \varphi_0&0&0&\hdots& 0\\
  0&\varphi_1&0& \hdots& 0\\
  %0&\varphi_1&\hdots&0\\
\vdots&\vdots&\vdots&\ddots&\vdots\\
  0&0&\hdots &\varphi_{m-2}&0}.
$$
This defines a morphism 
$\cM(G_0,\overline{\lieg}_1)\to \cM(G)$, whose image is in the fixed point locus of
the action of the group $\mu_m$ of $m$-th roots of unity on  $\cM(G)$
defined by sending $(E,\varphi)$ to $(E,\zeta\varphi)$, where $\zeta$ is a primitive  $m$-th root of unity. The Higgs bundle $(E,\varphi)$ is thus  a cyclic $\SL(n,\C)$-Higgs bundle as defined in Section \ref{vinberg-cyclic}.

The Hodge bundle defined by a quiver bundle (\ref{quiver-bundle}) with
$\varphi_{m-1}=0$ is known in this case as a
{\bf chain} (see \cite{AGS,GH,GHS}) and is a fixed point in $\cM(G)$ for the action of $\C^*$ given by rescaling the Higgs field.
The study of the Toledo invariant and Arakelov--Milnor inequality for chains and corresponding cyclic Higgs bundles as described in Section \ref{hodge-cyclic} is treated  in
detail in \cite{garcia-prada-gonzalez}.  We will just recall here the case  $m=2$, already studied in \cite{BGG03,BGG06}.

The  $m=2$ case corresponds to $V=V_0\oplus V_1$,
$$
G_0=\SSS(\GL(V_0)\times\GL(V_1)),\;\;\; \mbox{and}\;\;\;
\overline{\lieg}_1= \Hom(V_0,V_1)\oplus \Hom(V_1,V_0),
$$
and to the quiver representations described by the diagram

\[\begin{tikzcd} 	
	{V_0} &  & {V_1} 
		\arrow["{f_0}",curve={height=-12pt}, from=1-1, to=1-3]
	\arrow["{f_1}", curve={height=-12pt}, from=1-3, to=1-1].
	\end{tikzcd}\]

A
$(G^\theta,\lieg_1)$-Higgs pair over $X$ in the sense of Section \ref{higgs-pairs} is equivalent to a $4$-tuple
$(E_0,E_1,\varphi_0,\varphi_1)$, consisting of holomorphic vector bundles  $E_0$ and $E_1$ over $X$  of ranks 
$n_0=\dim  V_0$ and $n_1=\dim  V_1$, respectively with  $\det E_1= (\det E_0)^{-1}$, and homomorphisms
$$
\varphi_0: E_0\to E_1\otimes K_X\;\;\; \mbox{and}\;\;\; \varphi_1: E_1\to E_0\otimes K_X.
$$
In the terminology of \cite{AG}, this is equivalent to the $K_X$-twisted  quiver bundle
\[\begin{tikzcd} 	
	{E_0} &  & {E_1} 
		\arrow["{\varphi_0}",curve={height=-12pt}, from=1-1, to=1-3]
	\arrow["{\varphi_1}", curve={height=-12pt}, from=1-3, to=1-1],
	\end{tikzcd}\]
where the $K_X$-twisting means that the target of both homomorphisms are twisted by $K_X$.
Let
$(E_0,E_1,\varphi_0,\varphi_1)$ be  an element in $\cM(G_0,\overline{\lieg}_1)$, then we can associate to it the $\SL(n,\C)$-Higgs bundle  $(E,\varphi)$ with
$$
E=E_0\oplus E_1\;\;\; \mbox{and}\;\;\;
\varphi =
\begin{pmatrix}
  0 & \varphi_1 \\
  \varphi_0  & 0
\end{pmatrix}.
$$
This defines a morphism  $\cM(G^\theta,\lieg_1)\to \cM(G)$, whose image is in the fixed point locus of the involution of $\cM(G)\to \cM(G)$ defined by sending a $G$-Higgs bundle $(E,\varphi)$ to $(E,-\varphi)$ (see Section \ref{vinberg-cyclic}, \cite{garcia-prada-2007, garcia-prada-ramanan}).

The objects in  $\cM(G_0,\overline{\lieg}_1)$ are some times referred as $\SU(n_0,n_1)$-Higgs bundles since  by the non-abelian Hodge correspondence explained in
\ref{section-nahc}, the moduli space $\cM(G_0,\overline{\lieg}_1)$ is homemorphic to the $\SU(n_0,n_1)$-character variety of the fundamental group of $X$
(see \cite{BGG03,BGG06,garcia-prada-ramanan}). 
Indeed, for $G=\SL(n,\C)$ the classes in $\conj(G)/\sim$ corresponding to the trivial clique
are represented (see \cite{helgason} e.g.) by the conjugations
$\sigma_{n_0,n_1}$, with $0\leq n_0\leq n_1$ and $n_0+n_1=n$ given by 
$$
   \begin{aligned}
  \sigma_{n_0,n_1}:\SL(n,\C)   & \to \SL(n,\C) \\
    A  &\mapsto I_{n_0,n_1} (\overline{A}^t)^{-1} I_{n_0,n_1},
  \end{aligned}
$$
where
$$
I_{n_0,n_1} =
\begin{pmatrix}
  I_{n_0} & 0 \\
  0 & -I_{n_1}
\end{pmatrix}.
$$

One has $G^{\sigma_{n_0,n_1}}=\SU(n_0,n_1)$. In particular $\tau:=\sigma_{0,n}$ gives
the compact real form $\SU(n)$.
The elements in $\Aut_2(G)$ corresponding to $\sigma_{n_0,n_1}$ are given by
$\theta_{n_0,n_1}:=\tau\sigma_{n_0,n_1}$, and hence
$G^{\theta_{n_0,n_1}}=\SSS(\GL(n_0,\C)\times \GL(n_1,\C))$. Notice that the quasi-split real form in the trivial clique is $\SU(k,k)$ if $n=2k$, and $\SU(k,k+1)$ if $n=2k+1$.

A natural invariant associated to a $(G_0,\overline{\lieg}_1)$-Higgs pair $(E_0,E_1,\varphi_0,\varphi_1)$
is given by
$$
d=\deg E_0=-\deg E_1.
$$
The Toledo invariant in this case is $\tau=2d$ and,
if $(E_0,E_1,\varphi_0,\varphi_1)$ is semistable, one has
as a particular case of (\ref{mw-inequality})
$$
-\rank(\varphi_1) (g-1) \leq d \leq \rank(\varphi_0) (g-1),
$$
originally proved in  \cite{BGG03}, 
which implies the classical Milnor--Wood inequality
$$
 |d|\leq \min\{n_0,n_1\}(g-1).
$$

Let $\cM_d\subset \cM(G_0,\overline{\lieg}_1)$ the subvariety consisting of elements with fixed Toledo invariant $2d$. In \cite{BGG03,BGG06,BGGH} it is proved that
$\cM_d$ is connected if $|d|<\min\{n_0,n_1\}(g-1)$ and has
$2^{2g}$ connected components if $|d|=\min\{n_0,n_1\}(g-1)$. This result was proved by Hitchin for $\SU(1,1)$ in his original paper \cite{hitchin1987}.

A very special rigidity phenomenon takes place when $|d|=\min\{n_0,n_1\}(g-1)$. If $n_0\neq n_1$ this manifests in the fact that the dimension of $\cM_d$ is smaller than the expected dimension which is half of the dimension of $\cM(G)$
(see \cite{BGG03,BGG06} for details). When $n_0=n_1=m$  and
$|d|=m(g-1)$, then $\cM_d$ is isomorphic as an algebraic variety to the moduli space of $K_X^2$-twisted $(\hat{G},\liegl(m,\C))$-Higgs pairs over $X$ where
$$
\hat{G}=\{A\in \GL(m,\C)\;\;\mbox{such that}\;\; (\det A)^2=1\}.
$$

This isomorphism is known as {\bf Cayley correspondence} and takes place more generally for Hermitian real forms of {\bf tube type}
(see  \cite{BGG03,BGG06,BGR} for details). The $2^{2g}$ connected components of
$\cM_d$ in this case are known as {\bf Cayley components} and via the non-abelian Hodge correspondence are homeomorphic to components of the corresponding character variety consisting entirely of discrete and faithful representations. These are examples of {\bf higher Teichm\"uller spaces} (see \cite{wienhard,garcia-prada-handbook} for surveys on this subject and \cite{BCGGO,guichard-labourie-wienhard} for more recent results). In the case of $\SU(1,1)$, as described by Hitchin \cite{hitchin1987}, these components are identified with the usual Teichm\"uller space of the underlying smooth surface to $X$.
A similar rigidity phenomenon takes place also when the Toledo invariant is maximal  for Hodge bundles \cite{BCGT} and cyclic Higgs bundles \cite{garcia-prada-gonzalez}.

%%%%%%%%%%%%%%%%%%%%%%%%%%%%%%%%%%%%%%%%%%%%%%%%%%%%%%%%%%%%%%%%%%%%%%%%%%%%%%%%%%%%
\section{Vinberg $\theta$-pairs and the Hitchin fibration}\label{hitchin-fibration}
%%%%%%%%%%%%%%%%%%%%%%%%%%%%%%%%%%%%%%%%%%%%%%%%%%%%%%%%%%%%%%%%%%%%%%%%%%%%%%%%%%%%

Let $G$ be a semisimple complex Lie group and let $\theta\in \Aut(G)$ be of order $m$.
Consider the  $\Z/m\Z$-grading defined by $\theta$
   $$
   \lieg=\bigoplus_{i\in \Z/m\Z} \lieg_i,
   $$
   and the Vinberg $\theta$-pair $(G^\theta,\lieg_1)$.
Let $X$ be a compact Riemann surface and  $\cM(G^\theta,\lieg_1)$ be the moduli space 
of $(G^\theta,\lieg_1)$-Higgs pairs over $X$.
Recall from Sections \ref{vinberg-theta-pairs} and \ref{extended} (whose definitions and notions we follow)
that by \cite{vinberg}
   \begin{equation}\label{vinberg-theorem-2}
   \C[\lieg_1]^{G^\theta} \to \C[\liea]^{W(\liea)}=\C[f_1,\cdots,f_r],  \end{equation}
   where $r=\dim \liea=\rk(G^\theta,\lieg_1)$. Let $d_i=\deg f_i$. 

    Evaluating the polynomials $f_i$ on the Higgs field we have the  {\bf Hitchin map}
  \begin{equation}\label{hitchin-map}
    h: \cM(G^\theta,\lieg_1)\to B(G^\theta,\lieg_1)\cong \bigoplus_{i=1}^r
    H^0(X,K_X^{d_i}).
   \end{equation}

When $m=1$, (\ref{vinberg-theorem-2}) is the classical Chevalley restriction Theorem for the adjoint action of $G$ on $\lieg$, and 
one has  the original Hitchin map defined  in \cite{hitchin:duke} and the {\bf Hitchin integrable system}
  $$
  \cM(G)\to B(G)\cong \bigoplus_{i=1}^r
  H^0(X,K^{d_i}).
  $$
  Here $r=\rk G$ and  $\{d_1,\cdots, d_r\}$ are the {\bf exponents} of $G$.
  In \cite{hitchin:duke} Hitchin gives a description of the generic fibres for the  classical groups in terms of the Jacobian or
  appropriate Prym variety of a {\bf spectral curve}. 
  In \cite{donagi-gaitsgory} Donagi--Gaitsgory give a cameral curve description for  general $G$ as a {\bf gerbe} with
  generic abelian fibres. This approach was reformulated by  Ng\^o \cite{ngo}  in his proof of the
    {\bf Fundamental Lemma}.
    In  \cite{hitchin1992}  Hitchin constructed a section of the Hitchin map which he  identified with a connected component of the
    character variety  of the fundamental group of $X$  in a  split real form of $G$. This is called a {\bf Hitchin component}
    and is another instance of the higher Teichm\"uller spaces mentioned in Section \ref{hodge-cyclic}. 

    The case $m=2$ corresponds to symmetric pairs. Here, as explained in Section \ref{section-nahc},  the moduli space
    $\cM(G^\theta,\lieg_1)$
    of $(G^\theta,\lieg_1)$-Higgs pairs over $X$ is homeomorphic to the character variety of $\pi_1(X)$ in the real form $G^\sigma$,
    where $\sigma=\tau\theta$, and $\tau$ is a compact antiholomorphic involution of $G$ commuting with $\theta$. 
The Hitchin map in this situation has been studied by  Schaposnik for classical real forms using the  spectral curve approach, and 
   by   Pe\'on-Nieto \cite{peon} in terms of cameral covers for an arbitrary involution $\theta$.
    From both points of view, one can see that the generic fibres are abelian if and only if the
    real form is  quasi-split. This  non-abelian phenomenon is  very nicely illustrated for certain real forms in \cite{hitchin-schaposnik}. 
    A construction of the gerbe in the quasi-split case,  following the Donagi--Gaitsgory approach has been given in \cite{garcia-prada-peon}.
    This has also been studied  by  Leslie \cite{leslie} in relation to a  symmetric pair version of the Fundamental Lemma.
    Further study of the gerbe structure is being carried out by Hameister--Morrissey \cite{hameister-morrissey}.  
    A section of the Hitchin map in this case was constructed in \cite{garcia-prada-peon-ramanan} building upon the Kostant--Rallis section \cite{kostant-rallis} and the Hitchin--Kostant section given by Hitchin in \cite{hitchin1992}

 The study of the Hitchin map for general Vinberg $\theta$-pairs with $m>2$  is one of the main themes
    of \cite{garcia-prada-gonzalez}. There are several important  questions to address here. One is regarding the
    quasi-split condition for the pair $(G^\theta,\lieg_1)$ that was given in Section \ref{quasi-split-gradings}. While we have many examples of Vinberg $\theta$-pairs satisfying this condition, we are not aware of a classification of quasi-split pairs in this
    generality. Similarly to the case of symmetric pairs, this condition  ensures that the generic fibres are abelian. Another important issue  regards the construction of a section of the Hitchin map. This is  based on the
    existence of the Kostant--Weierstrass section for the quotient map $\lieg_1\to \lieg_1\sslash G^\theta$ mentioned in Section \ref{vinberg-theta-pairs}.

\providecommand{\bysame}{\leavevmode\hbox to3em{\hrulefill}\thinspace}
    
\bibliographystyle{amsplain}

\end{document}